\documentclass[12pt]{article}
\usepackage{amsmath,amssymb,amsthm, comment}

\newtheorem{theorem}{Theorem}[section]
\newtheorem{thmy}{Theorem}

\newtheorem{lemma}[theorem]{Lemma}
\newtheorem{corollary}[theorem]{Corollary}

\newtheorem{proposition}[theorem]{Proposition}
\newtheorem{example}[theorem]{Example}

\def\barr{\begin{array}}
\def\earr{\end{array}}

\title{Group cosets with all elements automorphically equivalent}
\author{Ryan McCulloch and Marius T\u arn\u auceanu}
\date{November 8, 2026}

\begin{document}

\maketitle

\begin{abstract}
Let $G$ be a finite group and let $N$ be a proper non-trivial normal subgroup of $G$. We say that $(G,N)$ is an \textit{automorphically equivalent pair} if for every element $x\in G\setminus N$, the elements in the coset $xN$ belong to the same orbit under the action of ${\rm Aut}(G)$. This generalizes the concept of a Camina pair and constitutes a particular type of equal order pair. In this note, we determine all automorphically equivalent pairs in the class of finite abelian groups.
\end{abstract}

{\small
\noindent
{\bf MSC2020\,:} Primary 20E34; Secondary 20D45, 20E45, 20K01.

\noindent
{\bf Key words\,:} automorphically equivalent pair, Camina pair, equal order pair, abelian group.}

\section{Introduction}

Let $G$ be a finite group and $N$ be a proper non-trivial normal subgroup of $G$. We say that $(G,N)$ is a \textit{Camina pair} if for every element $x\in G\setminus N$, the elements in the coset $xN$ are conjugate with $x$. Camina pairs have been introduced in \cite{3} and then extensively studied by many researchers. It is well-known that if $(G,N)$ is a Camina pair, then either $G$ is a Frobenius group with kernel $N$ (which happens if and only if the orders of $N$ and $G/N$ are coprime) or at least one among $N$ and $G/N$ is a group of prime-power order. 

A generalization of Camina pairs was recently introduced in \cite{4}. We say that $(G,N)$ is an \textit{equal order pair} if for every element $x\in G\setminus N$, the elements in the coset $xN$ have the same order as $x$. Obviously, Camina pairs are equal order pairs. Note that equal order pairs have many properties of Camina pairs, but with some important differences. They have been completely described for several particular classes of finite groups, such as finite abelian $p$-groups (see item (b) in Proposition 3.2 of \cite{4}).

\begin{thmy}
Let $G$ be a finite abelian group and $N$ be a proper non-trivial subgroup of $G$. Then $(G,N)$ is an equal order pair if and only if there exists a positive
integer $i$ such that $\Omega_{i-1}(G)\subseteq N\subseteq\Omega_i(G)$.
\end{thmy}

In what follows, we introduce a new concept, which is related to the above two. Given a finite group $G$ and a proper non-trivial normal subgroup $N$ of $G$, we call $(G,N)$ an \textit{automorphically equivalent pair} if for every element $x\in G\setminus N$, the elements in the coset $xN$ are automorphically equivalent, i.e. 
\begin{equation}
\forall\, x\in G\setminus N, \forall\, n\in N, \exists\, f\in{\rm Aut}(G) \mbox{ such that } xn=f(x).\nonumber
\end{equation}It is clear that a Camina pair is an automorphically equivalent pair and that an automorphically equivalent pair is an equal order pair. Moreover, these inclusions are strict, as the next example shows.

\begin{example}
Let $G=\mathbb{Z}_2\times\mathbb{Z}_8$, $N_1=\Omega_1(G)\cong\mathbb{Z}_2 \times \mathbb{Z}_2$ and $N_2=\Omega_2(G)\cong\mathbb{Z}_2\times\mathbb{Z}_4$. Then $(G,N_1)$ is an equal order pair, but not an automorphically equivalent pair, while $(G,N_2)$ is an automorphically equivalent pair, but not a Camina pair.
\end{example}

As a first step in investigating these pairs $(G,N)$, we focus here on finite abelian groups. In this setting, the problem easily reduces to $p$-groups. Our main result is as follows.

\begin{theorem}
Let $G\cong\mathbb{Z}_{p^{m_1}}\times\cdots\times\mathbb{Z}_{p^{m_k}}$ be a finite abelian $p$-group with $1\leq m_1\leq\dots\leq m_k$ and $N$ be a proper non-trivial subgroup of $G$ satisfying $\Omega_{i-1}(G)\subseteq N \subset \Omega_i(G)$ for some $i=1,\dots,m_k$.  Let $j_0$ be the smallest $j$ so that $m_j = m_k$.
\begin{itemize}
\item If $\Omega_i(\mathbb{Z}_{p^{m_{j_0}}}\times\cdots\times\mathbb{Z}_{p^{m_k}}) \subseteq N$, then $(G,N)$ is an automorphically equivalent pair if and only if $N \subseteq p^{m_k-i-1}G$.
\item If $\Omega_i(\mathbb{Z}_{p^{m_{j_0}}}\times\cdots\times\mathbb{Z}_{p^{m_k}}) \not\subseteq N$, then $(G,N)$ is an automorphically equivalent pair if and only if $N \subseteq p^{m_k-i}G$.
\end{itemize}
\end{theorem}

We present two examples illustrating Theorem 1.2.  Both examples were cross-verified in GAP \cite{8}.

\begin{example}
Let $G=\mathbb{Z}_2\times\mathbb{Z}_4\times\mathbb{Z}_8$.  Let $N = \Omega_1(\mathbb{Z}_2) \times \Omega_1(\mathbb{Z}_4)  \times \Omega_2(\mathbb{Z}_8)  \cong \mathbb{Z}_2\times\mathbb{Z}_2\times\mathbb{Z}_4$.  Then $\Omega_1(G) \subseteq N \subset \Omega_2(G)$.  So $i=2$, $m_k = 3$, and $\Omega_2(\mathbb{Z}_8) \subseteq N$.  We have $m_k - i - 1 = 0$, and $N \subseteq 2^0G = G$, thus $(G,N)$ is an automorphically equivalent pair.
\end{example}

\begin{example}
Let $G=\mathbb{Z}_8\times\mathbb{Z}_{32}$.  Let $N = \Omega_1(\mathbb{Z}_8) \times \Omega_1(\mathbb{Z}_{32})  \cong \mathbb{Z}_2\times\mathbb{Z}_2$.  Then $\Omega_1(G) \subseteq N \subset \Omega_2(G)$.  So $i=2$, $m_k = 5$, and $\Omega_2(\mathbb{Z}_{32}) \not\subseteq N$.  We have $m_k - i = 3$, and $N \not\subseteq 8G$, thus $(G,N)$ is not an automorphically equivalent pair.
\end{example}

By taking $i=m_k$ in Theorem 1.2, we have $\Omega_i(\mathbb{Z}_{p^{m_{j_0}}}\times\cdots\times\mathbb{Z}_{p^{m_k}}) \not\subseteq N$ as $N$ is a proper subgroup of $G$, and the condition $N \subseteq p^{m_k-i}G$ becomes $N \subseteq p^0G=G$, which is trivially true. Thus, the following corollary holds.

\begin{corollary}
Let $G$ be a finite abelian $p$-group of exponent $p^{m_k}$ and $N$ be a proper non-trivial subgroup of $G$ satisfying $\Omega_{m_k-1}(G)\subseteq N$. Then $(G,N)$ is an automorphically equivalent pair.
\end{corollary}

Finally, we propose a natural question regarding the above concept.

\bigskip\noindent{\bf Open problem.} Find automorphically equivalent pairs for other important classes of finite groups.
\bigskip

Most of our notation is standard and will usually not be repeated here. Elementary concepts and results on group theory can be found in \cite{6}.

\section{Background on isomorphic pairs of abelian groups}

We begin with a result of K. B\'uz\'asi \cite{2}.  See also \cite{7}.

\begin{thmy}
If $G$ is a finite abelian group and $H_1,H_2 \leq G$ are cyclic subgroups, then there exists $\varphi\in{\rm Aut}(G)$ such that $\varphi(H_1)=H_2$ if and only if $H_1 \cong H_2$ and $G/H_1\cong G/H_2$.   
\end{thmy}

In \cite{1} we find a strengthening of B\'uz\'asi's result to elements.

\begin{thmy}
Let $G$ be a finite abelian group and $x,y\in G$. Then there exists $\varphi\in{\rm Aut}(G)$ such that $\varphi(x)=y$ if and only
if $G/\langle x\rangle\cong G/\langle y\rangle$.
\end{thmy}

The following result appears in a paper of Hill \cite{5}.

\begin{thmy}
Suppose that $A$ and $B$ are subgroups of the finite abelian $p$-group $G$. The following conditions are equivalent.
\begin{itemize}
\item There exists $\varphi\in{\rm Aut}(G)$ such that $\varphi(A)=B$.
\item There is an isomorphism between the subgroups which preserves values (heights).
\end{itemize}
\end{thmy}

The $p$-height $h_p(x)$ of an element $x$ in a finite abelian $p$-group $G$ is the largest integer $r\geq 0$ such that $x\in p^rG$.  We say that the $p$-height of $0$ is $\infty$.  What Hill means by ``an isomorphism between the subgroups which preserves values (heights)'' is that there exists an isomorphism $f : A \rightarrow B$ so that $h_p(x) = h_p(f(x))$ for all $x \in A$.  

In the case when $A = \langle x \rangle$ and $B = \langle y \rangle$, $A \cong B$ if and only if $o(x) = o(y)$.  Furthermore, if $y'$ is another generator for $B$, then $y' = ay$ for some integer $a$ coprime to $p$, and we have $h_p(y) = h_p(y')$.  Hence if there exists an isomorphism $f : A \rightarrow B$ so that $h_p(x) = h_p(f(x))$ for all $x \in A$, then since $y' = f(x)$ is a generator for $B$, it follows that every isomorphism from $A$ to $B$ will be height preserving, so in particular the isomorphism $g : A \rightarrow B$ induced by $g(x)=y$ will be height preserving.  Finally, observe that if $o(x) = o(y)$ and $h_p(x) = h_p(y)$, then it follows that $h_p(p^nx) = h_p(p^ny)$ for all $n \geq 0$.  Thus, there exists a height preserving isomorphism from $A = \langle x \rangle$ to $B = \langle y \rangle$ if and only if $o(x) = o(y)$ and $h_p(x) = h_p(y)$.

Hence we obtain the following result that we need for the proof of Theorem 1.2.

\begin{lemma}
Let $G$ be a finite abelian $p$-group and $x,y\in G$. The following are equivalent.
\begin{itemize}
\item There exists $\varphi\in{\rm Aut}(G)$ such that $\varphi(x)=y$.
\item $G/\langle x\rangle\cong G/\langle y\rangle$.
\item We have $o(x)=o(y)$ and $h_p(x)=h_p(y)$.
\end{itemize}
\end{lemma}

\section{Proofs of the main results}

We begin with a proposition that will lead to a proof of Theorem 2.1.  Note that for $u,v$ in an abelian $p$-group $G$, we have $h_p(u) \neq h_p(v)$ implies $h_p(u+v)=\min\{h_p(u), h_p(v)\}$, and we have $h_p(u) = h_p(v)$ implies $h_p(u+v) \geq h_p(u)$.

\begin{proposition}
Let $G$ be a finite abelian $p$-group and $N$ be a proper non-trivial subgroup of $G$ satisfying $\Omega_{i-1}(G)\subseteq N\subseteq\Omega_i(G)$ for some positive integer $i$, and let $m = \mathrm{max}\{h_p(x) \,\, | \,\, x \in G \setminus N\}$. Then $(G,N)$ is an automorphically equivalent pair if and only if $N \subseteq p^{m}G$.
\end{proposition}

\begin{proof}
$(\Rightarrow)$ Assume $(G,N)$ is an automorphically equivalent pair. For every $x\in G \setminus N$ and $n\in N$, there exists an automorphism mapping $x$ to $x+n$, which by Lemma 2.1 requires $h_p(x)=h_p(x+n)$. If there were an element $n \in N$ such that $h_p(n) < h_p(x)$ for some $x \in G \setminus N$, then $h_p(x+n) = h_p(n) < h_p(x)$, a contradiction. Hence, we have
\begin{equation}\label{eq:height_bound}
    h_p(n)\geq h_p(x) \quad \text{for all } x \in G \setminus N \text{ and } n\in N. \nonumber
\end{equation}
Since $m = \mathrm{max}\{h_p(x) \,\, | \,\, x \in G \setminus N\}$, it follows that $N\subseteq p^{m}G$.

\medskip
\noindent $(\Leftarrow)$ Conversely, suppose that $N \subseteq p^{m} G$, and let $x\in G \setminus N$ and $n\in N$.  Since $\Omega_{i-1}(G)\subseteq N\subseteq\Omega_i(G)$ for some positive integer $i$, Theorem A implies that $o(x)=o(x+n)$.  Now $N \subseteq p^{m} G$ implies that $h_p(n)\geq m$ and $m = \mathrm{max}\{h_p(x) \,\, | \,\, x \in G \setminus N\}$ implies that $h_p(x) \le m$. 

Since $h_p(n) \geq h_p(x)$, we either have $h_p(n) > h_p(x)$, which immediately implies $h_p(x+n) = h_p(x)$, or $h_p(n) = h_p(x) = m$.  In the latter case, since $m = \mathrm{max}\{h_p(x) \,\, | \,\, x \in G \setminus N\}$, and since $x+n \notin N$, $h_p(x+n)$ cannot be greater than $h_p(x)$, and thus $h_p(x) = h_p(x+n)$.  And so Lemma 2.1 implies that $x$ and $x+n$ are automorphically equivalent.
\end{proof}

We are now able to prove our main result.

\begin{proof}[Proof of Theorem 1.2]
Let $m = \mathrm{max}\{h_p(x) \,\, | \,\, x \in G \setminus N\}$.  Suppose first that $\Omega_i(\mathbb{Z}_{p^{m_{j_0}}}\times\cdots\times\mathbb{Z}_{p^{m_k}}) \subseteq N$.  We work to show that $m = m_k - i - 1$.  We have $i < m_k$ since $N$ is a proper subgroup of $G$.  Let $x$ be an element of order $p^{i+1}$ in $\mathbb{Z}_{p^{m_k}}$.  Note that $x \notin N$ since $N \subset \Omega_i(G)$.  We have $h_p(x) = m_k - (i+1) =m_k - i - 1$.  We claim that this is the maximum.  Given $x = (x_1,\dots,x_k) \in G$, observe that $h_p(x) = \mathrm{min}\{ h_p(x_j) \,\, | \,\, j = 1,\dots,k \}$.  Let $x = (x_1,\dots,x_k) \in G \setminus N$ be arbitrary.  If $j < j_0$, then $h_p(x_j) \leq (m_k - 1) - i$ since $\Omega_{i-1}(G) \subseteq N$ and so $o(x_j) \geq p^i$.  If $j \geq j_0$, then $h_p(x_j) \leq m_k - (i+1)$ since $\Omega_i(\mathbb{Z}_{p^{m_{j_0}}}\times\cdots\times\mathbb{Z}_{p^{m_k}}) \subseteq N$ and so $o(x_j) \geq p^{i+1}$.  Thus $h_p(x) \leq m_k - i - 1$ and we have $m = m_k - i - 1$. 

Suppose now that $\Omega_i(\mathbb{Z}_{p^{m_{j_0}}}\times\cdots\times\mathbb{Z}_{p^{m_k}}) \not\subseteq N$.  We work to show that $m = m_k - i$.  In this case we have $x \in \mathbb{Z}_{p^{m_k}}$ of order $p^i$ so that $x \notin N$.  We have $h_p(x) = m_k - i$.  We claim that this is the maximum.   Let $x = (x_1,\dots,x_k) \in G \setminus N$ be arbitrary.  For any $j$ we have $h_p(x_j) \leq m_k  - i$ since $\Omega_{i-1}(G) \subseteq N$ and so $o(x_j) \geq p^i$.  Thus $h_p(x) \leq m_k - i$ and we have $m = m_k - i$.  

Theorem 1.2 then follows immediately from Proposition 3.1.
\end{proof}

\bigskip\noindent{\bf Funding.} The authors did not receive support from any organization for the submitted work. 

\bigskip\noindent{\bf Conflict of interest.} The authors declare that they have no conflict of interest. 

\bigskip\noindent{\bf Data availability.} The manuscript has no associated data.

\vspace*{3ex}
\small

\begin{minipage}[t]{8cm}
Ryan McCulloch\\
Department of Mathematics and Statistics\\
Binghamton University\\
Binghamton, NY 13902\\
e-mail: {\tt rmccullo1985@gmail.com}
\end{minipage}
\hfill
\begin{minipage}[t]{5cm}
Marius T\u arn\u auceanu \\
Faculty of  Mathematics \\
``Al.I. Cuza'' University \\
Ia\c si, Romania \\
e-mail: {\tt tarnauc@uaic.ro}
\end{minipage}

\end{document}